\documentclass[12pt]{article}

\usepackage{amsfonts}
\usepackage{amssymb}
\usepackage{amsmath}
\usepackage{enumerate}
\usepackage{stmaryrd}

\newtheorem{ttt}{Theorem}[section]
\newtheorem{llll}[ttt]{Lemma}
\newtheorem{ccc}[ttt]{Claim}
\newtheorem{eee}[ttt]{Example}
\newtheorem{rrr}[ttt]{Remark}
\newtheorem{sss}[ttt]{Statement}
\newtheorem{ddd}[ttt]{Definition}
\newtheorem{qqq}[ttt]{Question}
\newtheorem{cccc}[ttt]{Corollary}
\newtheorem{mmm}[ttt]{Main Lemma}
\newtheorem{ppp}[ttt]{Problem}

\newcommand{\bt}{\begin{ttt}}
\newcommand{\bl}{\begin{llll}}
\newcommand{\bc}{\begin{ccc}}
\newcommand{\bex}{\begin{eee}}
\newcommand{\br}{\begin{rrr}\upshape}
\newcommand{\bs}{\begin{sss}}
\newcommand{\bd}{\begin{ddd}\upshape}
\newcommand{\bq}{\begin{qqq}}
\newcommand{\bcor}{\begin{cccc}}
\newcommand{\bml}{\begin{mmm}}
\newcommand{\bpr}{\begin{ppp}}

\newcommand{\bp}{\noindent\textbf{Proof. }}

\newcommand{\et}{\end{ttt}}
\newcommand{\el}{\end{llll}}
\newcommand{\ec}{\end{ccc}}
\newcommand{\eex}{\end{eee}}
\newcommand{\er}{\end{rrr}}
\newcommand{\es}{\end{sss}}
\newcommand{\ed}{\end{ddd}}
\newcommand{\eq}{\end{qqq}}
\newcommand{\ecor}{\end{cccc}}
\newcommand{\eml}{\end{mmm}}
\newcommand{\epr}{\end{ppp}}

\newcommand{\qed}{\hspace{\stretch{1}}$\square$\medskip}

\newcommand{\lab}[1]{\label{#1}}

\numberwithin{equation}{section}

\newcommand{\imp}{\Rightarrow}
\newcommand{\emb}{\hookrightarrow}
\newcommand{\nemb}{\hookrightarrow\hspace{-11pt}\Arrownot\hspace{11pt}}
\renewcommand{\succ}{\textrm{succ}}
\newcommand{\dom}{\textrm{dom}}

\renewcommand{\es}{\emptyset}

\newcommand{\al}{\alpha}
\newcommand{\om}{\omega}
\newcommand{\si}{\sigma}
\newcommand{\ga}{\gamma}

\newcommand{\Baire}{\iB_1(\RR)}
\newcommand{\rest}{\mathord{\upharpoonright}} 

\newcommand{\AND}{\text{ and }}

\newcommand{\SetOf}[2]{\left\{#1 \ \left| \ #2 \right.\right\}}

\newcommand{\Sacks}{{\mathbb S}}
\newcommand{\Reals}{{\mathbb R}}
\newcommand{\Rationals}{{\mathbb Q}}
\newcommand{\Poset}{{\mathbb P}}
\newcommand{\Cantor}{{\mathbb C}}

\newcommand{\Naturals}{{\mathbb N}}

\newcommand{\NN}{\mathbb{N}}

\newcommand{\PP}{\mathbb{P}}
\newcommand{\QQ}{\mathbb{Q}}
\newcommand{\RR}{\mathbb{R}}
\renewcommand{\SS}{\mathbb{S}}
\newcommand{\CC}{\mathbb{C}}

\newcommand{\iB}{\mathcal{B}}

\newcommand{\iK}{\mathcal{K}}

\newcommand{\iP}{\mathcal{P}}

\begin{document}

\title{Chains of Baire class 1 functions and various notions of
special trees}

\author{M\'arton Elekes\thanks{Partially
supported by
Hungarian Scientific Foundation grant no.~37758 and F 43620.} 
\ and Juris Stepr\={a}ns\thanks{The second author's
research for this paper was partially supported by NSERC of Canada.}}

\maketitle

\begin{abstract}
Following Laczkovich we consider the partially ordered set
$\iB_1(\RR)$ of Baire class 1 functions endowed with the pointwise
order, and investigate the order types of the linearly ordered
subsets.  Answering a question of Komj\'ath and Kunen we show (in
$ZFC$) that special Aronszajn lines are embeddable into $\iB_1(\RR)$.
We also show that under Martin's Axiom a linearly ordered set
$\mathbb{L}$ with $|\mathbb{L}|<2^\omega$ is embeddable into
$\iB_1(\RR)$ iff $\mathbb{L}$ does not contain a copy of $\omega_1$ or
$\omega_1^*$.  We present a $ZFC$-example of a linear order of size
$2^\omega$ showing that this characterisation is not valid for orders
of size continuum.

These results are obtained using the notion of a compact-special tree;
that is, a tree that is embeddable into the class of compact subsets of
the reals partially ordered under reverse inclusion. 
We investigate how this notion is related to the well-known notion of
an $\RR$-special tree and also to some other notions of specialness.
\end{abstract}

\insert\footins{\footnotesize{MSC codes: Primary 26A21, 03E04; Secondary 
03E50, 03E15}}
\insert\footins{\footnotesize{Key Words: Baire class 1, pointwise ordering, 
linear ordering, Martin's Axiom, Aronszajn line, special tree}}

\section*{Introduction}

\bd
Given two partial orders $(\Poset,\leq_{\Poset})$ and 
$(\Poset',\leq_{\Poset'})$
the order $\Poset$  will be said to \textbf{embed into} 
$\Poset'$, denoted by 
$\Poset \emb \Poset'$, if there
is a mapping $\varphi \colon \Poset \to \Poset'$ such that 
$p_0 <_{\Poset} p_1$ implies $\varphi(p_0) <_{\Poset'} \varphi(p_1)$.
\ed

Note that this $\varphi$ need not be one-to-one in general, but for a 
linear order $\mathbb{L}$ the relation $\mathbb{L} \emb \Poset$ implies that
there is an order-isomorphic copy of $\mathbb{L}$ in $\PP$. As it is usual 
for trees, instead of $\Poset \emb \Poset'$ we will sometimes say that 
$\Poset$ is $\Poset'$-special.
From now on we will often write $\PP$ instead of $(\Poset,\leq_{\Poset})$ 
when there is no danger of confusion.

$\iB_1(\RR)$ is the class of Baire class 1 functions from $\RR$ to $\RR$; 
that is, pointwise limits of sequences of continuous real functions. 
This class
is partially ordered under the usual pointwise ordering; that is,
$f\leq g$ iff $f(r)\leq g(r)$ for every $r \in \RR$. Note that 
$f < g$ iff $f\leq g$ and $f(r) \neq g(r)$ for some $r \in \RR$.
The following problem was posed by Laczkovich.

\bpr\lab{p:Laczk}
Characterise those linear orders $\mathbb{L}$ for which $\mathbb{L} \emb 
\iB_1(\RR)$ holds.
\epr

What makes the Baire class 1 case particularly 
interesting is that the corresponding
questions for all other Baire classes are solved.
In the Baire class 0; that is, continuous case it is easy to see that 
$\mathbb{L} \emb \iB_0(\RR)$ iff $\mathbb{L} \emb \RR$,
while for $\al\geq 2$ Komj\'ath \cite{Ko} showed that even the question whether
$\om_2 \emb \iB_\al(\RR)$ is independent of $ZFC$.

Another motivation for Problem \ref{p:Laczk} may be that for 
the first sight it seems to be closely
related to the well-known theory of Rosenthal compacta. 
However, no direct connection has been found yet.

The earliest result concerning Problem \ref{p:Laczk} is a classical 
theorem of 
Kuratowski \cite[24.III.2']{Kur} stating that $\om_1 \nemb \iB_1(\RR)$. 
Note that $\al \emb \RR \emb \iB_1(\RR)$ for $\al < \om_1$. For some 
related results see \cite{EK}. 
It is shown in \cite{El} that, loosely speaking, starting from a class
of simple linear orders, say the finite ones, and applying all sorts of
countable operations one always obtains $\iB_1(\RR)$-embeddable linear 
orders. Therefore it is quite natural to guess that Kuratowski's theorem 
is the only restriction; that is, $\mathbb{L} \emb \iB_1(\RR)$ iff 
$\om_1,\om_1^* \nemb \mathbb{L}$. (Here $\om_1^*$ is the reversed $\om_1$.) 
However, Komj\'ath \cite{Ko} gave a consistent 
counterexample by showing in $ZFC$ that if $\mathbb{L}$ is a Souslin 
line then $\mathbb{L} \nemb \iB_1(\RR)$. But this still leaves open the 
possibility that the above answer to Laczkovich's problem is consistent
with $ZFC$.

\bq\lab{q:char}
Is it consistent that a linear order $\mathbb{L} \emb \iB_1(\RR)$ iff 
$\om_1,\om_1^* \nemb \mathbb{L}$?
\eq

Komj\'ath and Kunen independently asked the following natural question.

\bq\lab{q:Aron}
Is there an Aronszajn line $\mathbb{A}$ such that $\mathbb{A} \emb 
\iB_1(\RR)$?
\eq

In this paper we answer Question \ref{q:char} and Question \ref{q:Aron}. 

First we establish our basic tool in Section \ref{s:ML}, then make some
preparations in Section \ref{s:trees} by proving that nine notions of 
specialness coincide for countably branching trees. Then we answer Question 
\ref{q:Aron} in the positive in Section \ref{s:con}. More precisely, we
show that special Aronszajn lines are $\iB_1(\RR)$-embeddable, hence
there exists (in $ZFC$) a $\iB_1(\RR)$-embeddable Aronszajn line, and
consistently all Aronszajn lines are $\iB_1(\RR)$-embeddable. We also show 
in this section that under Martin's Axiom the characterisation in 
Question \ref{q:char} is valid for linear orders of cardinality strictly 
less than the continuum. In Section \ref{s:char} we answer Question 
\ref{q:char} in the negative (in $ZFC$). Finally, in Section \ref{s:open} 
we formulate some open problems.
 
The set-theoretic terminology followed in this paper can be found
e.g.~in \cite{Je} and \cite{Kun}. For an element $t$ of a tree
$\mathbb{T}$ denote $\succ(t)$ the set of immediate successors of
$t$. We say that a tree $\mathbb{T}$ is countably branching, if
$|\succ(t)| \leq \om$ for every $t \in \mathbb{T}$. All trees in this
paper are considered to be normal; that is, for $t_0, t_1 \in
\mathbb{T}$ the equation $\{t \in \mathbb{T} \colon t <_\mathbb{T} t_0\} =
\{t \in \mathbb{T} \colon t <_\mathbb{T} t_1\}$ implies $t_0=t_1$. The
basic facts about Baire class 1 functions can be found e.g.~in
\cite{Ke} or \cite{Kur}. An $F_\si$ set is a set that is the union of
countably many closed sets, a $G_\delta$ set is a set that is the
intersection of countably many open sets. 

\section{The main lemma}\lab{s:ML}

For a linear order $\mathbb{L}$, we say that $\mathbb{T}_\mathbb{L}$ 
is a (binary) partition tree of $\mathbb{L}$ (see \cite{To}), if it 
is constructed as follows.  
Denote by $\mathbb{T}_\al$ the $\al^\textrm{th}$ level of a tree
$\mathbb{T}$.  Elements of the partition tree will be nonempty
intervals; that is, convex subsets of $\mathbb{L}$, and the ordering
will be reverse inclusion.  Set
$(\mathbb{T}_\mathbb{L})_0=\{\mathbb{L}\}$. Once
$(\mathbb{T}_\mathbb{L})_\al$ is given, split every $I\in
(\mathbb{T}_\mathbb{L})_\al$ of at least two elements into two
disjoint nonempty intervals $I^+_0$ and $I^+_1$, and put
$(\mathbb{T}_\mathbb{L})_{\al+1} = \{I^+_i \colon I \in
(\mathbb{T}_\mathbb{L})_\al,\ |I|\geq 2, \ i\in 2\}$. We tacitly
assume that $I^+_0$ is the `left' interval; that is, for every $l_0
\in I^+_0$ and $l_1 \in I^+_1$ we have $l_0 \leq_\mathbb{L} l_1$.  For
$\al$ limit put $(\mathbb{T}_\mathbb{L})_\al = \{\cap_{\beta<\al}
I_\beta \colon I_\beta \in (\mathbb{T}_\mathbb{L})_\beta, \
\cap_{\beta<\al} I_\beta \neq \es \}$.

Denote by $\iK(\RR)$ the set of compact subsets of $\RR$ ordered 
under reverse inclusion.

\bd
We say that $\mathbb{T} \emb \iK(\RR)$ \textbf{strongly}, if there
exists an embedding which maps incomparable elements to disjoint sets;
that is, there exists an embedding $\varphi \colon \mathbb{T} \to \iK(\RR)$
such that $\varphi(t_0) \cap \varphi(t_1) = \es$ for every
$t\in\mathbb{T}$ and distinct $t_0, t_1 \in \succ(t)$.
\ed

\bml\lab{ML}
Let $\mathbb{L}$ be a linear order and $\mathbb{T}_\mathbb{L}$ be 
a partition tree of $\mathbb{L}$ such that
$\mathbb{T}_\mathbb{L}\emb\iK(\RR)$ strongly. Then
$\mathbb{L}\emb\iB_1(\RR)$.
\eml

\bp
Let $\varphi \colon \mathbb{T}_\mathbb{L} \to \iK(\RR)$ be a strong embedding. 
For every $l\in \mathbb{L}$ define
\[
A^l = \bigcup\{\varphi(I^+_0) \colon I \in \mathbb{T}_\mathbb{L}, \ |I|\geq 2,\ 
l \in I^+_1 \}.
\]
We claim that $\psi \colon \mathbb{L} \to \iB_1(\RR)$
\[
\psi(l) = \chi_{A^l}
\]
is the required embedding, where $\chi_H$ is the characteristic
function of the set $H$. As $\chi_{H_0} < \chi_{H_1}$ iff $H_0
\subsetneq H_1$, we first have to show that for $l_0 <_\mathbb{L} l_1$
the strict inclusion $A^{l_0} \subsetneq A^{l_1}$ holds.

Fix $l_0 <_\mathbb{L} l_1$.  First we show $A^{l_0} \subseteq
A^{l_1}$.  Suppose $I \in (\mathbb{T}_\mathbb{L})_\al, \ |I|\geq 2$
and $l_0 \in I^+_1$. We have to show that $\varphi(I^+_0) \subseteq
A^{l_1}$. There is a first level where $l_0$ and $l_1$ are not in the same
element of $\mathbb{T}_\mathbb{L}$, moreover, this is necessarily a 
successor level,
say $l_0, l_1 \in I^*\in (\mathbb{T}_\mathbb{L})_\al*$, $l_0 \in
(I^*)^+_0$ and $l_1 \in (I^*)^+_1$. Clearly, $\varphi((I^*)^+_0)
\subseteq A^{l_1}$.
If $\al < \al^*$ then $l_0 \in I^+_1$ implies $l_1 \in I^+_1$, hence
$\varphi(I^+_0) \subseteq A^{l_1}$. If $\al \geq \al^*$ then $I
\subseteq I^*$, hence $\varphi(I^+_0) \subseteq \varphi((I^*)^+_0)
\subseteq A^{l_1}$.

Now we show $A^{l_0} \neq A^{l_1}$. By compactness, $C =
\bigcap\{\varphi(I) \colon l_0 \in I \in \mathbb{T}_\mathbb{L}\} \neq \es$. Using
$\varphi((I^*)^+_0) \subseteq A^{l_1}$ again, we obtain $C \subseteq
A^{l_1}$. We claim that $C \cap A^{l_0} = \es$. In order to show this
we have to check that $l_0 \in I^+_1$ implies $\varphi(I^+_0) \cap C =
\es$. But this is clear, as $C \subseteq \varphi(I^+_1)$ and $\varphi$ 
is a strong embedding.

What remains to be shown is that $\chi_{A^l} \in \iB_1(\RR)$ for every
$l \in \mathbb{L}$. A characteristic function $\chi_H$ is of Baire
class 1 iff $H$ is simultaneously $F_\sigma$ and $G_\delta$, hence we
have to check this for $A^l$.  It is well known (see \cite[22.27]{Ke}
or \cite[24.III.1]{Kur}) that if for some $\xi<\om_1$ the nonincreasing
transfinite sequences $\{F_\al\}_{\al<\xi}$ and $\{H_\al\}_{\al<\xi}$
of closed subsets of $\RR$ satisfy $F_\al \supseteq H_\al$ for every
$\al<\xi$ and $H_\al \supseteq F_\beta$ for every $\al<\beta<\xi$,
then the set
\[
\bigcup_{\al<\xi} (F_\al\setminus H_\al)
\]
is simultaneously $F_\sigma$ and $G_\delta$.

Fix $l \in \mathbb{L}$. Let $\xi^l$ be the ordinal for which $\{l\}
\in (\mathbb{T}_\mathbb{L})_{\xi^l}$ holds. As every strictly
decreasing transfinite sequence of compact subsets of $\RR$ is
countable, $\xi^l < \om_1$. For $\al < \xi^l$ the unique interval
$I\in(\mathbb{T}_\mathbb{L})_\al$ with $l \in I$ has at least two
elements, so define
\[
F^l_{\al+1} = H^l_{\al+1} = \varphi(I^+_0) \cup \varphi(I^+_1),
\]
if $l \in I^+_0$, and
\[
F^l_{\al+1} = \varphi(I^+_0) \cup \varphi(I^+_1),
\]
\[
H^l_{\al+1} = \varphi(I^+_1)
\]
if $l \in I^+_1$.  For $\al < \xi^l$ limit, which includes the case
$\al=0$, define
\[
F^l_{\al}=H^l_{\al}=\varphi(I).
\]
Clearly, $F^l_\al \supseteq H^l_\al$ for every $\al<\xi^l$ and it is
easy to see that $F^l_\al \supseteq F^l_\beta$ and $H^l_\al \supseteq
H^l_\beta$ for every $\al<\beta<\xi^l$. Using that $F^l_\al$ is 
monotone nonincreasing, in order to obtain that $H^l_\al \supseteq 
F^l_\beta$ for every $\al<\beta<\xi^l$ it is sufficient to check that 
$H^l_\al \supseteq F^l_{\al+1}$ for every $\al<\xi^l$, which is 
straightforward. Therefore $\bigcup_{\al<\xi^l} (F^l_\al\setminus
H^l_\al)$ is $F_\sigma$ and $G_\delta$. Using that our embedding 
$\varphi$ is strong we obtain
\[
A^l = \bigcup_{\al<\xi^l} (F^l_\al\setminus H^l_\al),
\]
so the proof is complete.
\qed

\section{Various notions of special trees}\lab{s:trees}

In this section we prove that the relation $\mathbb{T}_\mathbb{L} \emb 
\iK(\RR)$ strongly can be translated to $\mathbb{T}_\mathbb{L} \emb 
\RR$. As specialness of trees is interesting in its own right, we 
prove that, at least for countably branching trees, this is also 
equivalent to specialness in certain other senses.
Let $\Cantor$ denote the Cantor set (not the complex plane!) 
with its inherited ordering as a
subset of $\Reals$. The Prikry-Silver partial order will be denoted by
$\Sacks$ -- it consists of all partial functions $f\colon\Naturals \to 2 =
\{0,1\}$ with co-infinite domain ordered under inclusion. 

\bd
We say that $\mathbb{T} \emb \SS$ \textbf{strongly}, if there
exists an embedding which maps incomparable elements to incompatible 
functions; that is, there exists an embedding $\varphi \colon \mathbb{T} 
\to \SS$ such that for every $t\in\mathbb{T}$ and distinct 
$t_0, t_1 \in \succ(t)$ 
there exists $n \in \dom(\varphi(t_0)) \cap \dom(\varphi(t_1))$ such that 
$\varphi(t_0)(n) \neq \varphi(t_1)(n)$. 
\ed

\bt\lab{specials}
Let $\mathbb{T}$ be a countably branching tree, e.g.~a partition tree. Then
the following are equivalent. 
\begin{enumerate}[(1)]
    \item\label{C} $\mathbb{T}$ is $\Cantor$-special (Cantor-special)
    \item\label{R} $\mathbb{T}$ is $\Reals$-special
    \item\label{sS} $\mathbb{T}$ is strongly $\Sacks$-embeddable
    \item\label{sKC} $\mathbb{T}$ is strongly $\iK(\Cantor)$-embeddable
    \item\label{KC} $\mathbb{T}$ is $\iK(\Cantor)$-special
    \item\label{sKR} $\mathbb{T}$ is strongly $\iK(\RR)$-embeddable
    \item\label{KR} $\mathbb{T}$ is $\iK(\RR)$-special
    \item\label{PN} $\mathbb{T}$ is $(\iP(\Naturals),\subseteq)$-special
    \item\label{S} $\mathbb{T}$ is $\Sacks$-special        
    \end{enumerate}
\et

\bp
(\ref{C}) $\imp$ (\ref{R}): This is immediate.

(\ref{R}) $\imp$ (\ref{sS}): Let $\varphi\colon \mathbb{T} \to \RR$ be an
embedding, and let $\{q_n\colon n\in\NN\}$ enumerate $\QQ$.  Set
$\dom(\psi(t)) = \{n\in\NN\colon  q_n < \varphi(t)\}$ and define $\psi(t)\colon 
\{n\in\NN\colon  q_n < \varphi(t)\} \to 2$ by induction along $\mathbb{T}$
as follows.  At limit nodes simply let $\psi(t)$ be the union of all
$\psi(s)$ such that $s \subset t$. Given that $\psi(t)$ is defined,
enumerate $\succ(t)$ as $\{t_k\colon  k\in\NN\}$, and by induction on $k$
pick distinct $n_k\in\NN$ such that $\varphi(t) \leq q_{n_k} <
\varphi(t_k)$. For $n\in \NN$ such that $q_n < \varphi(t_k)$ set
$\psi(t_k)(n) = \psi(t)(n)$ if $q_n < \varphi(t)$, $\psi(t_k)(n_k) =
1$, and $\psi(t_k)(n) = 0$ otherwise. It is easy to check that $\psi\colon 
\mathbb{T}\to \SS$ is a strong embedding.

(\ref{sS}) $\imp$ (\ref{sKC}):  Let $\varphi\colon  \mathbb{T} \to \SS$ be a
strong embedding. Identify $\CC$ with $2^\NN$; that is, the set of
functions from $\NN$ to $2$. For $t\in \mathbb{T}$ define $\psi(t) =
\{f\in 2^\NN \colon  \varphi(t) \subseteq f\}$. Then $\psi\colon \mathbb{T}\to
\iK(\Cantor)$ is a strong embedding.

(\ref{sKC}) $\imp$ (\ref{KC}): Obvious.

(\ref{KC}) $\imp$ (\ref{C}): Again, identify $\CC$ with $2^\NN$. Let
$\{g_n\}_{n=1}^\infty$ enumerate all $g\colon k\to 2$ where $k\in
\Naturals$ and send $K\in \iK(\Cantor)$ to 
\[
\sum_{\substack{n\in\NN \\ \not{\exists} f\in K \ g_n\subseteq f}} 
\frac{2}{3^{n+1}}.
\]

(\ref{sKC}) $\imp$ (\ref{sKR}): Obvious.

(\ref{sKR}) $\imp$ (\ref{KR}): Obvious.

(\ref{KR}) $\imp$ (\ref{R}): Enumerate $\{(p,q) \colon  p,q\in\QQ, \ p<q\}$ as
$\{(p_n,q_n)\colon n\in\NN\}$, and send $K\subseteq \RR$ to $\sum_{(p_n,q_n) \cap
K = \es} \frac{1}{2^n}$.

(\ref{R}) $\imp$ (\ref{PN}): Enumerate $\QQ$ as $\{q_n\colon n\in\NN\}$, and send 
$r\in\RR$ to $\{n\in\NN: q_n < r\}$. 

(\ref{PN}) $\imp$ (\ref{S}):  Send $H\subseteq\Naturals$ to the function 
that is constant $0$ on $\{2n\colon n\in H\}$ and undefined elsewhere.

(\ref{S}) $\imp$ (\ref{R}): Send $f\in\SS$ to 
$-\sum_{n\notin\dom(f)} \frac{1}{2^n}$.
\qed

\br The assumption that the tree $\mathbb{T}$ is countably branching
cannot be dropped, as if $\succ(t)$ has cardinality larger than the
continuum for some $t\in \mathbb{T}$ then $\mathbb{T}$ is clearly not
strongly $\iK(\Cantor)$-special but it can be $\RR$-special.

It is well-known, that even for countably branching trees 
$\QQ$-specialness is not equivalent to the properties listed in the above
theorem. Indeed, one can show that $\sigma \QQ$ (see the proof of 
Theorem \ref{nochar} or \cite{To}) is $\RR$-special, but not $\QQ$-special.

It is also well-known, that for $\omega_1$-trees (trees of height $\omega_1$
with countable levels) it is independent of $ZFC$ whether $\RR$-specialness is
equivalent to $\QQ$-specialness. Indeed, for one direction it is enough that 
under $MA$ all $\omega_1$-trees with no uncountable branches are
$\QQ$-special, which was shown e.g.~in \cite{BM}. The other direction was
proved by Baumgartner (see e.g.~\cite{Sh}), who constructed an $\RR$-special,
non-$\QQ$-special Aronszajn tree under $\diamondsuit$.
\er

\section{Consequences for $\iB_1$-embeddability}\lab{s:con}

In this section we answer Question \ref{q:Aron} and give an affirmative 
answer to Question \ref{q:char} in the case $|\mathbb{L}|<2^\om$.

\bt
Let $\mathbb{A}$ be a special Aronszajn line; that is, for some partition tree
$\mathbb{T}_\mathbb{A}$ of $\mathbb{A}$ we have
$\mathbb{T}_\mathbb{A} \emb \QQ$. Then $\mathbb{A} \emb \iB_1(\RR)$.
\et

\bp
Clearly, $\mathbb{T}_\mathbb{A} \emb \RR$, hence Theorem
\ref{specials} yields $\mathbb{T}_\mathbb{A} \emb \iK(\RR)$ strongly,
therefore by the Main Lemma \ref{ML} we obtain $\mathbb{A} \emb
\iB_1(\RR)$.
\qed

\bt
Assume Martin's Axiom. Then for a linear order $\mathbb{L}$ with
$|\mathbb{L}|<2^\om$ the relation $\mathbb{L} \emb \iB_1(\RR)$ holds
iff $\om_1,\om_1^* \nemb \mathbb{L}$.
\et

\bp
First suppose $\om_1 \emb \mathbb{L}$ or $\om_1^* \emb \mathbb{L}$. By
the theorem of Kuratowski \cite[24.III.2']{Kur} every strictly monotone
transfinite sequence in $\iB_1(\RR)$ is countable, hence $\mathbb{L}
\nemb \iB_1(\RR)$.  Now suppose $\om_1,\om_1^* \nemb \mathbb{L}$. 
It
follows that there is no strictly decreasing sequence of subintervals
of $\mathbb{L}$ of length $\om_1$, hence $\mathbb{T}_\mathbb{L}$ has
at most $\om_1$ levels, where of course $\mathbb{T}_\mathbb{L}$ is 
a partition tree of $\mathbb{L}$. Each level of this tree 
is a disjoint family of nonempty
intervals of $\mathbb{L}$, so $|\mathbb{L}|<2^\om$ implies
$|(\mathbb{T}_\mathbb{L})_\al| < 2^\om$ for every $\al$.  By Martin's
Axiom $\om_1 < 2^\om$ and $2^\om$ is regular, therefore
$|\mathbb{T}_\mathbb{L}| < 2^\om$.  Under Martin's Axiom every tree of
cardinality less than $2^\om$ with no branch of length $\om_1$ is
$\QQ$-special \cite{BM}, hence $\mathbb{T}_\mathbb{L} \emb \QQ$, and
we can repeat the previous proof.
\qed

\section{Answer to Question \ref{q:char}}\lab{s:char}

Now we answer Question \ref{q:char} in the negative, using some ideas from
\cite{To}.

\bt\lab{nochar}
There exists a linear order $\mathbb{L}$ such that $\om_1,\om_1^*
\nemb \mathbb{L}$ but still $\mathbb{L} \nemb \iB_1(\RR)$.
\et

\bp
Define
\[
\si\Baire = \{l \colon  \xi<\om_1,\ l\colon \xi \to \Baire \textrm{ strictly 
increasing}\}.
\]
This set becomes a tree if we partially order it by extension; that is, 
$l_0 \leq_\mathbb{T} l_1$ iff $l_0 \subseteq l_1$.

\bl
$(\si\Baire,\leq_\mathbb{T}) \nemb \Baire$.
\el

\bp
Suppose $\varphi \colon  \si\Baire \to \Baire$ is an embedding. Then the 
transfinite recursion
\[
l^*(\al) = \varphi(l^*\rest\al)
\]
produces a strictly increasing sequence of length $\om_1$ in $\Baire$, 
which is impossible by Kuratowski's theorem \cite[24.III.2']{Kur}.
\qed

This lemma shows that in order to finish the proof of Theorem
\ref{nochar} it is sufficient to construct a linear order
$\leq_\mathbb{L}$ on $\si\Baire$ extending $\leq_\mathbb{T}$ such that
$\om_1, \om_1^* \nemb (\si\Baire,\leq_\mathbb{L})$. So fix an
arbitrary bijection $\Phi \colon  \Baire \to \RR$ and define
$\leq_\mathbb{L}$ to be the usual lexicographical ordering as follows. 
The functions $l_0 \colon  \xi^{l_0} \to \Baire$ and $l_1 \colon  \xi^{l_1} \to
\Baire$ are incomparable with respect to $\leq_\mathbb{T}$ iff there
exists $\al < \xi^{l_0}, \xi^{l_1}$ such that $l_0(\al) \neq
l_1(\al)$. In such a case choose the minimal such $\al$ and define
$l_0 <_\mathbb{L} l_1$ iff $\Phi(l_0(\al)) < \Phi(l_1(\al))$.

Now we prove that $\om_1, \om_1^* \nemb (\si\Baire,\leq_\mathbb{L})$. 
Suppose $\{l_\eta\}_{\eta<\om_1}$ is strictly monotonic. 
We prove by induction on $\beta<\om_1$ that there exists $l^* \colon  \om_1 
\to \Baire$ such that 
for every $\beta<\om_1$ there exists $\eta_\beta$ 
such that for $\eta \geq \eta_\beta$
\[
l_\eta(\beta) = l^*(\beta).
\]
Suppose this holds for every $\ga<\beta$. If $\eta \geq
\sup\{\eta_\gamma \colon  \ga<\beta\}$ then $l_\eta\rest\beta =
l^*\rest\beta$, and hence $\Phi(l_\eta(\beta))$ is monotonic in $\RR$,
and therefore is constant above some $\eta_\beta$. As $\Phi$ is a
bijection, $l_\eta(\beta)$ is also constant for $\eta \geq
\eta_\beta$. Defining $l^*(\beta) = l_{\eta_\beta}(\beta)$ finishes
the induction. But once again, the existence of the strictly monotone
sequence $\{l^*(\al)\}_{\al<\om_1}$ contradicts Kuratowski's theorem.
\qed

\section{Open questions}\lab{s:open}

The fundamental open problem is still of course Problem \ref{p:Laczk}. 
However, we formulate here a couple of related questions.

We mentioned in the Introduction that, starting from some simple
linear orders, countable operations always result in
$\Baire$-embeddable orders. However, we do not know whether the class
of $\Baire$-embeddable orders itself is closed under these
operations. It is shown in \cite{El} that the answer is affirmative
for all these operations provided that it is affirmative for the
simplest such operation, namely, for the operation that doubles the
points of the order. That is why we are particularly interested in the
following.

\bq
Suppose $\mathbb{L} \emb \Baire$, where $\mathbb{L}$ is a linear order. 
Does $\mathbb{L} \times \{0,1\} \emb \Baire$, where the ordering of 
$\mathbb{L} \times \{0,1\}$ is the usual lexicographical order?
\eq

Denote $\Delta^0_2(\RR)$ the class of subsets of $\RR$ that are 
simultaneously $F_\si$ and $G_\delta$. The ordering is reverse inclusion.
Clearly,
\[
\QQ \emb \RR \emb \Delta^0_2(\RR) \emb \Baire,
\]
and it can be shown that the first two arrows cannot be reversed.
How about the third one?

\bq
$\Baire \emb \Delta^0_2(\RR)$?
\eq

\bq
Suppose $\mathbb{L} \emb \Baire$, where $\mathbb{L}$ is a linear order. 
Does $\mathbb{L} \emb \Delta^0_2(\RR)$? How about trees instead of 
linear orderings?
\eq

\bigskip
\bigskip

\noindent
\textsc{R\'enyi Alfr\'ed Institute, Re\'altanoda u.~13-15, H-1053, Budapest, Hungary}

\textit{Email address}: \verb+emarci@renyi.hu+

URL: http://www.renyi.hu/\~{}emarci

\bigskip

\noindent
\textsc{Department of Mathematics, York University, 4700 Keele Street, 
Toronto, Ontario, Canada M3J 1P3}

\textit{Email address}: \verb+steprans@yorku.ca+

\end{document}